\documentclass[11pt,preprint]{article}

\usepackage{amsmath}
\usepackage{epsfig}
\usepackage{graphicx}
\usepackage{amssymb}
\usepackage{amsthm}
\usepackage{latexsym}
\usepackage{color}

 

\newcommand{\ra}{\rightarrow}

\newcommand{\E}{\mathsf{E}}
\newcommand{\R}{\mathbb{R}}
\newcommand{\Pe}{\mathbb{P}}
\newcommand{\rd}{\mathrm{d}}

\newcommand{\Lp}{{{L}_2}}
\newcommand{\Real}{\mathfrak{R}}
\newcommand{\Img}{\mathfrak{I}}
\newcommand{\C}{\mathbb{C}}

\newtheorem{theorem}{Theorem}[section]
\newtheorem{lemma}[theorem]{Lemma}
\newtheorem{proposition}[theorem]{Proposition}
\newtheorem{corollary}[theorem]{Corollary}

\newtheorem{definition}[theorem]{Definition}

\newtheorem{remark}[theorem]{Remark}

\title{A comparative linear mean-square stability analysis of Maruyama- and Milstein-type methods}

\author{Evelyn Buckwar and Thorsten Sickenberger\footnote{Support by the
Leverhulme Trust through project ``Stability issues for SDEs'' (F/00 276/K) is gratefully acknowledged.} \\[1em]
{\scriptsize Maxwell Institute and Heriot-Watt University, Dept.~of Mathematics, Edinburgh, EH14 4AS, Scotland, UK}\\[-0.3em]
{\scriptsize E-mail: \texttt{\{e.buckwar/t.sickenberger\}@hw.ac.uk}}\\
} 

\date{}

\begin{document}

\maketitle

\begin{abstract}
In this article we compare the mean-square stability properties of
the $\theta$-Maruyama and $\theta$-Milstein method that are used to solve stochastic
differential equations. For the linear stability analysis, we propose an extension of the standard geometric Brownian motion as a test equation and consider a scalar linear test equation with several multiplicative noise terms. This test equation allows to begin investigating the influence of multi-dimensional noise on the stability behaviour of the methods while the analysis is still tractable. Our findings include: (i) the stability condition for the $\theta$-Milstein method and thus, for some choices of $\theta$, the conditions on the step-size, are much more restrictive than those for the $\theta$-Maruyama method; (ii) the precise stability region of the $\theta$-Milstein method explicitly depends on the noise terms. Further, we investigate the effect of introducing  partially implicitness in the diffusion approximation terms of Milstein-type methods, thus obtaining the possibility to control the stability properties of these methods with a further method parameter $\sigma$. Numerical examples illustrate the results and provide a comparison of the stability behaviour of the different methods.



\end{abstract}

\section{Introduction}

In recent years the area of numerical analysis of stochastic differential equations (SDEs) has expanded at a fast pace. This
interest has been driven by different application areas, such as computational finance, neuroscience or electrical circuit engineering. A large part of research in stochastic numerics has been aimed towards the development and strong and weak convergence analysis of several classes of numerical methods. A further important issue for the investigation of numerical methods consists of examining methods for their ability to preserve qualitative features of the continuous system they are developed to approximate. A linear stability analysis is usually the first step of an analysis in this direction. For this the method of interest is applied to a scalar linear test equation and stability conditions on method parameters and step-size are derived and compared with the stability condition for the test equation. In deterministic numerical analysis the underlying idea for a linear stability analysis is based on the following line of reasoning: one linearises and centres a nonlinear ordinary differential equation $\mathbf{x}'(t)=f(t,\mathbf{x})$ around an equilibrium, the resulting linear system $\mathbf{x}'(t)= A \mathbf{x}(t)$ ($A$ the Jacobian of $f$ evaluated at equilibrium) is then diagonalised and the system thus decoupled, justifying the use of the scalar test equation $x'(t)=\lambda x(t),\,\lambda \in \C$, for the analysis. We refer to, for example,~\cite[Chapter IV.2]{HaWa} for more detail on this procedure.

\medskip%
In the stochastic case the same fundamental problem exists: we wish to preserve the qualitative behaviour of solutions of nonlinear stochastic differential equations following discretisation by a numerical method. Once again the starting point is a linear stability analysis. We can linearise and centre a nonlinear SDE around an equilibrium solution (see~\cite{Ha80} for the corresponding theory), and this procedure yields a system of SDEs with an $m$-dimensional driving Wiener process of the form ${\rd}X(t)= FX(t){\rd}t + \sum_{r=1}^m G_r X(t)\,{\rd}W_r(t), t> 0\,.$
Research on stability analysis for SDEs has focused on the scalar linear test equation ${\rd}X(t)= \lambda X(t){\rd}t + \mu_1 X(t) {\rd}W_1(t)$, where its solution is called geometric Brownian motion. This corresponds to considering the linearised SDE system when it can be completely decoupled and with one driving Wiener process. Relevant references are given by, e.g.,~\cite{BHW06,DH00,DH00a,SaMi96}. Some first explorations of a linear stability analysis for systems of SDEs have been made in~\cite{RathBal,SaMi02}, and in~\cite{BK10} suitable linear test systems have been derived.

\medskip%
The most common and well-known methods to treat SDEs numerically are the Euler-Maruyama approximation and the Milstein scheme developed in the last century. These methods and their drift-implicit counterparts have been studied for their stability behaviour, see for example~\cite{ArtAv,DH00,DH00a,SaMi96}, however always only for a one-dimensional noise term in the test equation. Our aim in this article is to investigate the influence of multi-dimensional noise on the mean-square stability behaviour of these numerical methods and to compare their stability behaviour. Assuming that the matrices $F$ and $G_r$ can be diagonalised and the system above decoupled, we propose a scalar linear equation with an $m$-dimensional Wiener process, that is $\rd X(t) = \lambda X(t) \rd t + \sum_{r=1}^m \mu_r X(t) \rd W_r(t)$, as a suitable test equation. We consider the $\theta$-Maruyama and the $\theta$-Milstein method and study the asymptotic mean-square stability properties of these methods. We find that the stability conditions for the $\theta$-Milstein method are not only stronger than those for the $\theta$-Maruyama method, but they also become more restrictive when increasing the number of noise terms in the SDE. In particular we find it impossible to conclude $A$-stability for the $\theta$-Milstein method, in the sense that we can not define a value or bound $\theta_{\text{bound}}$ for the parameter $\theta$, such that the method applied to the test equation would be $A$-stable for all $\theta\geq \theta_{\text{bound}}$ {\em for all $m\geq1$ and all parameters $\mu_r$, $r=1,\ldots,m$}. This indicates that analysing the stability behaviour of numerical methods using test equations with only a one-dimensional Wiener process may not provide sufficient insight into the properties of the methods for practical high-dimensional application problems. We then present a modification of the $\theta$-Milstein method, where we introduce partial implicitness involving a further parameter $\sigma$ in the hope to have a better control on the stability behaviour of the method. The stability analysis indicates that this is the case, although the stability condition still depends on the number of noise terms. Numerical experiments illustrate the influence that the choice of method and method parameters have for practical simulation tasks.

\medskip%
In Section~\ref{s.cc} we introduce the linear test equation, the $\theta$-Maruyama method and the $\theta$-Milstein scheme, we give a precise definition of asymptotic mean-square stability and state stability conditions for the test equation in the continuous case. Section~\ref{s.dc} is devoted to the analysis of the stability properties of the methods in terms of the arising stochastic difference equations and we compare the stability regions of both types of methods in Section~\ref{s.comp}. In Section~\ref{s.sigma} we introduce a new class of Milstein type methods and analyse their stability behaviour. We illustrate the theoretical findings by some numerical experiments in Section~\ref{s.num}.

\section{Preliminaries} \label{s.cc}

In this section we introduce the stochastic differential and difference equations, as well as the notions of stability, that we consider in this article.

\medskip%
We will be concerned with asymptotic {\em mean-square stability} of the zero
solution of a test equation with respect to perturbations in the
initial data. As a test equation we employ the following scalar
linear stochastic differential equation with multiplicative noise
\begin{equation}\label{testsde}
 \rd X(t) = \lambda X(t) \rd t + \sum_{r=1}^m \mu_r X(t) \rd W_r(t), \quad t \geq t_0 \geq 0, \quad X(t_0)= X_0,
\end{equation}
driven by an $m$-dimensional standard Wiener process $W(t) =
(W_1(t), \dots, W_m(t))$ given on the probability space
$(\Omega,\mathcal{F},\Pe)$ with a filtration
$(\mathcal{F}_t)_{t\geq t_0}$. By $\E$ we denote expectation with respect to $\Pe$. We assume that the coefficients
$\lambda$ and $\mu_r$, $(r=1,\dots, m)$ are complex-valued and
without loss of generality suppose that the initial value $X_0$ is non-random (\cite[Chapter 4]{Mao1}). The solution $(X(t))_{t\geq t_0}$ of~\eqref{testsde} is a complex-valued stochastic process, where $W$ is a vector-valued Wiener process in $\R^m$. Alternatively this complex-valued SDE can be written
as 2-dimensional vector SDE in the real and imaginary parts. We denote the real and imaginary part of a complex number $x$ by $ \Real(x)$ and $ \Img(x) $, respectively, $\bar{x}$ denotes the complex conjugate and $| x |$ stands for the absolute value of some $x \in \C$.

\medskip %
For any non-zero initial value
Equation~\eqref{testsde} has a non-trivial path-wise unique strong
solution which we denote by $X(t;t_0,X_0)$ when we wish to
emphasise its dependence on the initial data. For $X_0=0$ the
equation obviously admits the \emph{zero solution}, which is a steady
state solution of the equation.

\begin{definition} (cf. \cite{Ha80,Mao1}) \label{d.stabil}
The zero solution of Equation \eqref{testsde}, $X(t) \equiv 0$, is
\begin{enumerate}
\item
 {\em  mean-square stable}, if for each $\epsilon >0$, there exists a $\delta \geq 0$ such that
\[    \E (|X(t;t_0,X_0)|^2) < \epsilon  \]
whenever $t\geq t_0$ and $ |X_0| < \delta$;
\item {\em asymptotically mean-square stable}, if it is
mean-square stable and if there exists a $\delta \geq 0$ such that
whenever $|X_0| < \delta$
\[    \E (|X(t;t_0,X_0)|^2) \ra 0 \ \mbox{for}\ t \ra \infty\,. \]
\end{enumerate}
\end{definition}

The zero solution is called unstable if it is not stable in the mean-square sense. Definition~\ref{d.stabil} is slightly
more general than necessary in the present context, as for the simple linear equation given by~\eqref{testsde} we can take $\delta$ arbitrarily large and thus it does not play a significant role.

\begin{proposition} \label{l.contprob} \cite{Ha80,Mao1} The zero solution of Equation~\eqref{testsde}
is asymptotically mean-square stable
if and only if
\begin{equation}\label{zerostabcond}
 \Real(\lambda) + \frac12 \sum_{r=1}^m |\mu_r|^2 < 0\,.
\end{equation}
\end{proposition}

We now discuss the stochastic difference equations that arise by
applying the $\theta$-Maruyama method and the $\theta$-Milstein
scheme to the scalar test equation~\eqref{testsde}.
We consider numerical methods for computing approximations $X_i
\approx X(t_i)$ of the solution of the test equation~\eqref{testsde}
at discrete time points $t_i = ih$ ($i=0,1,\dots$) with constant step-size $h$.

\medskip%
The $\theta$-Maruyama method applied to the test equation~\eqref{testsde} reads
\begin{equation}\label{thetamethod}
X_{i+1} = X_i + h \left( \theta \lambda X_{i+1} + (1-\theta)
\lambda X_i \right) + \sqrt{h} \sum_{r=1}^m \mu_r X_i
\,\xi_{r,i}\,,\quad i=0,1,\ldots\,,
\end{equation}
where we have replaced the Wiener increments $I_r^{t_i,t_i+h} = W_r(t_i+h)-W_r(t_i)$ by the scaled random variables $\sqrt{h}\,
\xi_{r,i}$. Here each $\{\xi_{r,i}\}_{i\in\mathbb{N}}$ is one of $m$ independent sequences of mutually independent standard Gaussian random variables, i.e., each $\xi_{r,i}$ is ${\mathcal N}(0,1)$-distributed.

\medskip%
The $\theta$-Milstein method applied to the test equation~\eqref{testsde} is given by
\begin{eqnarray}
\nonumber X_{i+1} &=& X_i + h \left( \theta \lambda X_{i+1}
+ (1-\theta) \lambda X_i \right) + \sqrt{h} \sum_{r=1}^m \mu_r X_i\, \xi_{r,i}\\
\label{milsteinmethod}&& + \frac12\, h \sum_{r=1}^m \mu_r^2 X_i\,
(\xi_{r,i}^2 - 1) + \frac12\, h \!\!\!\sum_{r_1, r_2=1 \atop r_1 \neq r_2}^m \mu_{r_1} \mu_{r_2} X_i\,
\xi_{{r_1},i}\xi_{{r_2},i} \,,\quad i=0,1,\ldots\,.
\end{eqnarray}
Here we have additionally replaced the multiple Wiener integrals
$$I_{r_1,r_2}^{t_i,t_i+h}:=\int_{t_i}^{t_i+h} \int_{t_i}^{s}\rd  W_{r_2}(u)\, \rd W_{r_1}(s)$$ as follows: (i) If $r_1 = r_2$ the multiple Wiener integral $I_{r,r}^{t_i,t_i+h}$ can be replaced by $\frac12 h(\xi_{r,i}^2 - 1)$ for $r=1,\ldots,m$, and (ii) if $r_1 \neq r_2$ we use the identity $I_{r_1,r_2} + I_{r_2,r_1} = I_{r_1} I_{r_2} $ to obtain
$$\mu_{r_1} \mu_{r_2} I_{r_1,r_2} + \mu_{r_2} \mu_{r_1} I_{r_2,r_1} = \mu_{r_1} \mu_{r_2} ( I_{r_1,r_2} + I_{r_2,r_1}) = \mu_{r_1} \mu_{r_2}  I_{r_1}I_{r_2} = \mu_{r_1} \mu_{r_2} h \xi_{{r_1},i}\xi_{{r_2},i}\,. $$

\medskip %
For $\theta = 0$ the method~\eqref{thetamethod} reduces to the
\emph{(forward) Euler-Maruyama scheme},
which is explicit in the drift as well as in the diffusion part.
For $\theta > 0$ the methods~\eqref{thetamethod} and
\eqref{milsteinmethod} are drift-implicit. For $\theta = 1/2$ the
scheme~\eqref{thetamethod} is known as~\emph{stochastic
trapezoidal rule} and for $\theta = 1$ we obtain the
\emph{backward Euler-Maruyama method}.

\medskip %
The two last terms in the $\theta$-Milstein method~\eqref{milsteinmethod} represent a higher order approximation of the diffusion part in Equation~\eqref{testsde} in the sense of mean-square convergence. We call a method \emph{mean-square
convergent with order $\gamma$} ($\gamma>0$) if the global error,
$X(t_i) - X_i$, satisfies
\begin{equation*}
 \max_{i=0,1,\ldots} \| X(t_i) - X_i \|_\Lp := \max_{i=0,1,\ldots} (\E |X(t_i) - X_i|^2 )^{1/2} \leq C h^\gamma \ \text{ as } \ h \ra 0\,,
\end{equation*}
with a positive error constant $C$, which is independent of the step-size $h$.
It is well-known that in the case of multiplicative noise the
$\theta$-Milstein method \eqref{milsteinmethod} is mean-square convergent of order $\gamma=1$, whereas the $\theta$-Maruyama method \eqref{thetamethod} is mean-square convergent of order $\gamma=1/2$.

\medskip %
Obviously the stochastic difference equations~\eqref{thetamethod} and
\eqref{milsteinmethod} admit the zero solution for the initial
value $X_0=0$, which is a steady state solution as well. For any
non-zero initial value $X_0$, the equations have a unique solution
provided $(1-h\,\lambda\,\theta) \neq 0$. We write $X_i(t_0,X_0)$ again when we want to emphasise that the solution of the difference equations depends on the initial data.

\begin{definition} \label{d.stabildiff}
The zero solution of the difference equations~\eqref{thetamethod}
and~\eqref{milsteinmethod} is
\begin{enumerate}
\item
 {\em  mean-square stable}, if for each $\epsilon >0$, there exists a $\delta \geq 0$ such that
\[    \E (|X_i(t_0,X_0)|^2) < \epsilon  \]
whenever $i\geq 0$ and $|X_0| < \delta$;
\item {\em asymptotically mean-square stable}, if it is
mean-square stable and if there exists a $\delta \geq 0$ such that
whenever $|X_0| < \delta$
\begin{equation}\label{stabdefmethod}
 \E (|X_i(t_0,X_0)|^2) \ra 0 \quad \mbox{for}\quad i \ra \infty\,.
\end{equation}
\end{enumerate}
\end{definition}

Stability conditions will now involve the coefficients $\lambda$, $\mu_r$ $(r=1,\dots,m$) of the test
equation \eqref{testsde}, as well as the method parameter $\theta$ and the step-size $h$.

\medskip%
It will be useful to describe the stability region of a stochastic differential or difference equation.
We follow the presentation in \cite{DH00a} and consider sets of
parameters $S_\text{SDE}$, $S_{\theta\text{-Mar}}(\theta,h)$,
$S_{\theta\text{-Mil}}(\theta,h)$ for which the zero solutions of
the continuous and the discrete equations are asymptotically stable, that is
\begin{eqnarray*}
S_\text{SDE} &:=& \{ \lambda, \mu_1,\ldots,\mu_m \in \C:\
\mbox{\eqref{zerostabcond} holds} \}, \\
S_{\theta\text{-Mar}}(\theta,h)&:=& \{ \lambda,
\mu_1,\ldots,\mu_m \in \C :\
\mbox{\eqref{stabdefmethod} holds for solutions of \eqref{thetamethod}} \}, \\
S_{\theta\text{-Mil}}(\theta,h)&:=& \{ \lambda,
\mu_1,\ldots,\mu_m \in \C : \ \mbox{\eqref{stabdefmethod} holds for solutions of \eqref{milsteinmethod}} \}.
\end{eqnarray*}

Further, we consider the extension of the deterministic notion of A-stability~\cite{HaWa,Lamb} to the mean-square analysis setting and say that a numerical method is \emph{A-stable} in mean-square~\cite{DH00,DH00a}, if whenever the zero solution of
\eqref{testsde} is asymptotically mean-square stable, then the same is true for the zero
solution of the method  for any step-size $h> 0$. Using the above definitions of the stability regions, we call the $\theta$-Maruyama or $\theta$-Milstein method A-stable in mean-square if for all $h > 0$ we have $S_\text{SDE} \subseteq S_{\theta\text{-Mar}}(\theta,h)$ or $S_\text{SDE} \subseteq S_{\theta\text{-Mil}}(\theta,h)$.

\section{Stability analysis for the numerical methods}\label{s.dc}

For the stochastic difference equations~\eqref{thetamethod} and~\eqref{milsteinmethod} we now derive
stability conditions in dependence on the method parameter
$\theta$ and the applied step-size $h$, and compare these conditions with
those for the continuous problem given in Lemma~\ref{l.contprob}.

\subsection{The $\theta$-Maruyama method}

We start by rearranging the stochastic difference equation~\eqref{thetamethod} into the following one-step recurrence equation. Then by squaring and taking expectations we obtain a recurrence equation for the second moments $\E |X_i|^2$. As
mentioned before, we have to assume that $(1- \theta h \lambda) \neq 0$ to guarantee the existence of a unique solution of the
recurrence equations.

\medskip %
Rearranging Equation~\eqref{thetamethod} yields the recurrence
\begin{equation}\label{thetarecurrence}
X_{i+1} = \bigg(a + \sum_{r=1}^m b_r \,\xi_{r,i} \bigg)
X_i\,,\quad i=0,\,1,\ldots\,,
\end{equation}
where
\begin{equation}\label{thetaabc}
 a := 1 + \frac{h \lambda}{1 - \theta h \lambda}\qquad \mbox{and}  \qquad b_r := \frac{\sqrt{h} \mu_r}{1- \theta h
 \lambda}\,.
\end{equation}
Then the recurrence equation for the second moment $\E|X_i|^2$, using $\E(\xi_{r,i})
= 0$, $\E(\xi_{r,i}^2) = 1$, and the complex conjugate equation of \eqref{thetarecurrence}, reads
\begin{eqnarray}
 \E |X_{i+1}|^2 &=& \E \bigg( \big(a + \sum_{r=1}^m b_r \,\xi_{r,i}\big) \big(\bar a + \sum_{r=1}^m \bar b_r \,\xi_{r,i}\big) \bigg) \E  |X_i|^2 \nonumber\\
 &=& \bigg( |a|^2 + \sum_{r=1}^m |b_r|^2 \bigg) \E
 |X_i|^2\,. \label{sqthetarecurrence}
\end{eqnarray}

\medskip%
From~\eqref{sqthetarecurrence} we can immediately read off necessary and sufficient stability conditions in terms of the above parameters.

\begin{lemma} \label{lemMarstababc}
The zero solution of the recurrence equation~\eqref{thetarecurrence}, representing a one-step Maruyama-type method applied to the test equation~\eqref{testsde}, is asymptotically mean-square stable if and only if
\begin{equation}\label{thetaineqabc}
  |a|^2 + \sum_{r=1}^m |b_r|^2 < 1 \,.
\end{equation}
\end{lemma}
Now rewriting the stability conditions~\eqref{thetaineqabc} in terms of the parameters $\lambda$ and
$\mu_r$ $(r=1,\dots,m)$ in the test equation~\eqref{testsde}, the
method parameter $\theta$ and the step-size $h$, using~\eqref{thetaabc}, we obtain
the following result.
\begin{corollary}
The zero solution of the stochastic difference equation given by the $\theta$-Maruyama method~\eqref{thetamethod} applied to the scalar linear test equation~\eqref{testsde} is asymptotically
mean-square stable if and only if
\begin{eqnarray}\label{stabcondtheta}
 \Real(\lambda) + \frac12 \sum_{r=1}^m  |\mu_r|^2 + \frac12 h ( 1- 2\theta) |\lambda|^2  &<&
 0\,.
\end{eqnarray}
\end{corollary}
We note that the first two terms in the left-hand sides of~\eqref{stabcondtheta} are equal to
the left-hand side of~\eqref{zerostabcond}, that is they
correspond to the stability condition for the continuous problem.

\medskip%
Now comparing the stability condition for the
continuous problem to that of the discrete problem, we immediately
obtain from~\eqref{stabcondtheta} an extension of the result
\cite[Thm.~4.1]{DH00a} to the case of~\eqref{testsde} driven by a
multi-dimensional Wiener process.

\begin{corollary}\label{c.Maruyama}
For all $h> 0$ it holds that
\[ \begin{array}{rcll}
S_{\theta\emph{-Mar}}(\theta,h)&\subset & S_\emph{SDE}& \quad \mbox{for}\quad  0 \leq \theta < 1/2\,,\\
S_{\theta\emph{-Mar}}(\theta,h)&=& S_\emph{SDE} & \quad \mbox{for} \quad  \theta = 1/2\,,\\
S_{\theta\emph{-Mar}}(\theta,h) &\supset& S_\emph{SDE}& \quad
\mbox{for} \quad \theta > 1/2\,.
\end{array}\]
In particular, for $\,\theta \geq 1/2$ the $\theta$-Maruyama method
is A-stable in mean-square. For $\,0 \leq \theta < 1/2$ and $(\lambda,\mu_1, \dots, \mu_m) \in S_\emph{SDE}$, the stability condition \eqref{stabcondtheta} for the $\theta$-Maruyama method is satisfied if and only if
\begin{equation}
 h < \frac{ -2 (\Real(\lambda) + \frac12 \sum_{r=1}^m |\mu_r|^2)}{(1-2\theta)|\lambda|^2} \,.
\end{equation}
\end{corollary}

\subsection{The $\theta$-Milstein method}

We now turn to the $\theta$-Milstein method and first follow the same steps as in the previous section to deduce a recurrence equation for the second moments of the solution of~\eqref{milsteinmethod} and read off the corresponding stability conditions. Rearranging the difference equation \eqref{milsteinmethod}, we
obtain
\begin{equation}\label{milsteinrecurrence}
X_{i+1} = \bigg(\widehat a + \sum_{r=1}^m b_r \,\xi_{r,i} +
\sum_{r_1,r_2=1}^m c_{r_1,r_2} \,\xi_{r_1,i}\xi_{r_2,i}\bigg) X_i,
\end{equation}
\begin{equation}\label{milsteinabc}
\!\!\!\!\!\!\!\!\!\!\!\!\!\!\!\!\!\!\!\!\!\!\!\!\!\!\!\!\!\!\mbox{where} \qquad \qquad  \widehat a := a - \sum_{r=1}^m c_{r,r}, \quad b_r := \frac{\sqrt{h} \mu_r}{1- \theta h \lambda}, \quad c_{r_1,r_2} = \frac{\frac12 h \mu_{r_1} \mu_{r_2}}{1- \theta h \lambda},
\end{equation}
and $a$ is given in \eqref{thetaabc}. Note that rewriting the parameter $\widehat a$ in terms of $a$ provides a convenient way to compare the stability conditions for the $\theta$-Maruyama and $\theta$-Milstein method.

\medskip%
With analogous calculations as in the previous section and by additionally using  $\E(\xi_{r,i}^3) = 0$ and $\E(\xi_{r,i}^4) = 3$, we find the recurrence for the second moments of the $\theta$-Milstein method as
\begin{eqnarray}
\nonumber \lefteqn{\E |X_{i+1}|^2 }&&\\
\nonumber&=& \E \bigg( \!\Big( \widehat a + \sum_{r=1}^m b_r \,\xi_{r,i} +
\sum_{r_1,r_2=1}^m c_{r_1,r_2} \,\xi_{r_1,i}\xi_{r_2,i}\Big) \,\Big( \bar {\widehat a} + \sum_{r=1}^m \bar b_r \,\xi_{r,i} +
\sum_{r_1,r_2=1}^m \bar c_{r_1,r_2} \,\xi_{r_1,i}\xi_{r_2,i}\Big) \!\bigg) \E |X_i|^2\\
\nonumber &=& \bigg( \widehat a \bar {\widehat a} + \widehat a \sum_{r=1}^m \bar c_{r,r} +  \bar {\widehat a} \sum_{r=1}^m c_{r,r}
+ \sum_{r=1}^m b_r \bar b_r + \Big(\sum_{r=1}^m c_{r,r} \Big) \Big(\sum_{r=1}^m \bar c_{r,r} \Big) \\
\nonumber && \qquad + 2 \sum_{r=1}^m c_{r,r} \bar c_{r,r} + \! \sum_{r_1, r_2=1 \atop r_1 \neq r_2 }^m\!\! c_{r_1,r_2} \cdot \!\! \sum_{r_1, r_2=1 \atop r_1 \neq r_2 }^m\!\! \bar c_{r_1,r_2}   \bigg)  \E |X_i|^2\\
\nonumber &=& \bigg( \! \Big( \widehat a + \sum_{r=1}^m c_{r,r} \Big)\Big( \bar {\widehat a} + \sum_{r=1}^m \bar c_{r,r} \Big)
+ \sum_{r=1}^m b_r \bar b_r \, + 2 \sum_{r=1}^m c_{r,r} \bar c_{r,r} + \! \sum_{r_1, r_2=1 \atop r_1 \neq r_2 }^m\!\! c_{r_1,r_2} \cdot \!\! \sum_{r_1, r_2=1 \atop r_1 \neq r_2 }^m\!\! \bar c_{r_1,r_2}   \bigg)  \E |X_i|^2\\
\label{milsteinrecur} &=& \bigg( |a|^2 + \sum_{r=1}^m |b_r|^2 + 2
\sum_{r=1}^m |c_{r,r}|^2  + \! \sum_{r_1, r_2=1 \atop r_1 \neq r_2 }^m\!\! c_{r_1,r_2} \cdot \!\! \sum_{r_1, r_2=1 \atop r_1 \neq r_2 }^m\!\! \bar c_{r_1,r_2}   \bigg)  \E |X_i|^2\,.
\end{eqnarray}

Again we can read off from \eqref{milsteinrecur} necessary and sufficient stability conditions in terms of the parameters $a$, $b_r$ and $c_{r_1,r_2}$.

\begin{lemma} \label{lemMilstababc}
The zero solution of the recurrence equation \eqref{milsteinrecurrence} given by the $\theta$-Milstein method applied to the test equation \eqref{testsde}, is asymptotically mean-square stable if and only if
\begin{equation}\label{milsteinineqabc}
  |a|^2 + \sum_{r=1}^m |b_r|^2 + 2 \sum_{r=1}^m |c_{r,r}|^2 + | c_\text{sum} |^2 < 1\,,
\end{equation}
where $c_\text{sum} = \sum_{r_1, r_2=1,\, r_1 \neq r_2 }^m  c_{r_1,r_2}$.
\end{lemma}

\begin{remark}
When applied to the linear scalar stochastic differential equation~\eqref{testsde}, other variants
of one-step Maruyama-type methods and one-step Milstein-type methods can
quite often be rearranged into the recurrence equations~\eqref{thetarecurrence} and~\eqref{milsteinrecurrence},
respectively, with appropriate definitions of the parameters $a$,
$b_r$ and $c_{r_1,r_2}$. Then Lemmata~\ref{lemMarstababc} and~\ref{lemMilstababc} can be applied to obtain stability conditions interpreted with these definitions of $a$, $b_r$ and $c_{r_1,r_2}$.
\end{remark}

The next corollary follows by rewriting the stability condition~\eqref{milsteinineqabc} in terms of the original parameters  using~\eqref{milsteinabc}.
\begin{corollary} \label{CorMilStab}
The zero solution of the stochastic difference equation given by the $\theta$-Milstein method~\eqref{milsteinmethod} applied to the scalar linear test equation~\eqref{testsde} is asymptotically mean-square stable if and only if
\begin{eqnarray}\label{stabcondmilstein}
 \Real(\lambda) + \frac12 \sum_{r=1}^m  |\mu_r|^2 + \frac12 h ( 1- 2\theta) |\lambda|^2 + \frac14 h \sum_{r=1}^m |\mu_r|^4 + \frac18 \, h \,\, \Big| \!\!\! \sum_{r_1, r_2=1 \atop r_1 \neq r_2 }^m\!\! \mu_{r_1} \mu_{r_2} \Big|^2 &<&
 0\,.
\end{eqnarray}
\end{corollary}
Again the first two terms in the left-hand side of~\eqref{stabcondmilstein} are equal to
the left-hand side of~\eqref{zerostabcond}. However, when comparing the above condition with the stability condition~\eqref{stabcondtheta} for the $\theta$-Maruyama method, we see that the left hand side of~\eqref{stabcondmilstein} contains two additional terms, which are always non-negative and depend on the noise intensities $\mu_{r}$ ($r=1,\dots,m$) but are independent of the parameter $\theta$. Thus the precise stability region $S_{\theta\emph{-Mil}}(\theta,h)$ of the $\theta$-Milstein method depends on the noise intensities $\mu_r$ ($r=1,\dots,m$). Suppose our aim is to use the method parameter $\theta$ to determine the optimal stability region of the $\theta$-Milstein method for a given set of parameters in~\eqref{testsde}, that is we aim to find a value of $\theta$ such that the sum of the third to the fifth term in~\eqref{stabcondmilstein} vanish for any step-size $h$. Then, assuming that~\eqref{zerostabcond} holds, we obtain for every set of parameters $\mu_1, \ldots, \mu_m$ a different value for that optimal $\theta_{\mu_1, \ldots, \mu_m}$. We define
\begin{equation} \label{thetaopt}
\theta_{\mu_1, \ldots, \mu_m} \,:=\, \frac12 + \sum_{r=1}^{m}
\frac{|\mu_r|^4}{4\,|\lambda|^2} + \frac{|\mu_\text{sum}|^2}{8\, |\lambda|^2}\,,
\end{equation}
where
\begin{equation}\label{musum}
\mu_\text{sum} := \sum_{r_1, r_2=1,\, r_1 \neq r_2 }^m \mu_{r_1} \mu_{r_2}\,.
\end{equation}
For this optimal $\theta_{\mu_1, \ldots, \mu_m}$ we immediately obtain from the stability condition~\eqref{stabcondmilstein} for the $\theta$-Milstein
method the following corollary.
\begin{corollary}\label{c.Milstein}
For all $h> 0$ it holds that
\[ \begin{array}{rcll}
S_{\theta\emph{-Mil}}(\theta,h)&\subset & S_\text{SDE}& \quad \mbox{for} \quad  0 \leq \theta < \theta_{\mu_1, \ldots, \mu_m}\\
S_{\theta\emph{-Mil}}(\theta,h)&=& S_\text{SDE} & \quad \mbox{for} \quad \theta = \theta_{\mu_1, \ldots, \mu_m}\\
S_{\theta\emph{-Mil}}(\theta,h) &\supset& S_\text{SDE}& \quad
\mbox{for}\quad  \theta > \theta_{\mu_1, \ldots, \mu_m}\,.
\end{array}\]
Assuming that $(\lambda,\mu_1, \dots, \mu_m) \in S_\emph{SDE}$, for $\,0 \leq \theta <\theta_{\mu_1, \ldots, \mu_m}$, the stability condition \eqref{stabcondmilstein} for the $\theta$-Milstein method is satisfied if and only if
\begin{equation}\label{hrestrictionMil}
 h < \frac{ -2 (\Real(\lambda) + \frac12 \sum_{r=1}^m |\mu_r|^2)}{(1-2\theta)|\lambda|^2 + \frac12 \sum_{r=1}^m |\mu_r|^4 + \frac14 |\mu_\text{sum}|^2} \,.
\end{equation}
In particular, condition \eqref{hrestrictionMil} implies that for $\theta> \frac12$ and $ (1-2\theta)|\lambda|^2 + \frac12 \sum_{r=1}^m |\mu_r|^4 + \frac14 |\mu_\text{sum}|^2< 0 $, the stability condition \eqref{stabcondmilstein} for the $\theta$-Milstein method is satisfied for any step-size $h > 0$\,. \\
In the case that $(\lambda,\mu_1, \dots, \mu_m) \not\in S_\emph{SDE}$, and $\theta> \frac12$ and $ (1-2\theta)|\lambda|^2 + \frac12 \sum_{r=1}^m |\mu_r|^4 + \frac14 |\mu_\text{sum}|^2> 0 $, then the zero solution of the $\theta$-Milstein scheme is also unstable only if condition \eqref{hrestrictionMil} is imposed on the step-size $h$.
\end{corollary}

\medskip
Intrinsic in the concept of $A$-stability is the idea that the property of $A$-stability of a method holds for the whole class of differential equations considered as test equations. In the setting of this article this implies that we would need to find a $\theta_\text{bound}$, a value of or a bound on  $\theta_{\mu_1, \ldots, \mu_m}$, which is independent of $\lambda, \mu_1, \ldots, \mu_m$, such that $S_{\theta\emph{-Mil}}(\theta,h) \supseteq S_\text{SDE}$ for $ \theta > \theta_\text{bound}$. In the case of the $\theta$-Maruyama method the corresponding value of $\theta_\text{bound}$ is $\frac12$, see Corollary \ref{c.Maruyama}.

However, considering the simple case of multi-dimensional noise terms having equal noise intensities, we find that, at best, we can find an upper bound $\theta_\text{bound}$ independent of the given parameters $\lambda, \mu_1, \ldots, \mu_m$, but depending on the number of noise sources. To see this, let $\mu_1 = \mu_2 = \cdots = \mu_m$. Then, using the squared stability inequality \eqref{zerostabcond} in the form $\frac14 |\mu_1|^4 / |\lambda|^2 < \frac1{m^2}$ and $\sum\limits_{r_1, r_2=1,\, r_1 \neq r_2 }^m \mu_{1} \mu_{1}= \mu_{1}^2 (m^2-m)$ in \eqref{thetaopt} yields
\begin{eqnarray}
\nonumber \theta_{\mu_1, \ldots, \mu_m} &<& \frac12 + \frac1m + \frac18 \frac{|\mu_1|^4}{|\lambda|^2} (m^2-m)^2\\
&<& \frac12 + \frac1m + \frac12 (m-1)^2\,,
\end{eqnarray}
Thus for example, for one, two and three noise sources ($m=1,\,2,\,3$), we can find a $\theta_\text{bound}$ as $\frac32,\,\frac32$ and $\frac{15}{6}$, respectively. In the one-dimensional noise case this corresponds to the result in~\cite[Cor.~2.2]{DH00}. But in general there exists no upper bound $\theta_\text{bound}$ for $\theta_{\mu_1, \ldots, \mu_m}$ such that the $\theta$-Milstein method for any $\theta \geq \theta_\text{bound}$ is $A$-stable for the whole class of equations \eqref{testsde} with arbitrary many noise terms. In the general case of the SDE \eqref{testsde} it is possible, essentially using the Cauchy-Schwarz inequality several times, to derive an upper bound  $\theta_\text{bound}$ for $\theta_{\mu_1, \ldots, \mu_m}$ which also only depends on $m$, but as it involves $m^3$ it becomes pointless for large $m$. Thus, it appears that the stability condition~\eqref{stabcondmilstein} for the $\theta$-Milstein method becomes very restrictive for an increasing number of noise terms.

\section{Comparisons of the stability regions}\label{s.comp}

In this section we aim to illustrate the results of the previous sections by visually comparing the stability regions  $S_{\theta\emph{-Mar}}(\theta,h)$ and $S_{\theta\emph{-Mil}}(\theta,h)$. As there are too many parameters in the system to do so in a two-dimensional plot, we only consider the case of the test SDE~\eqref{testsde} with multi-dimensional noise where all the terms have the same noise intensity, i.e.,  $\mu_1 = \mu_2 = \dots = \mu_m$, and real-valued coefficients, i.e., we have $\lambda, \mu_1 \in \R$. We essentially follow~\cite{DH00a} regarding the scaling of the parameters in the plots. Thus we set
\[
 x:=h \lambda \quad \text{and} \quad y := h m \mu_1^2 \,.
\]
The stability conditions~\eqref{zerostabcond}, \eqref{stabcondtheta} and~\eqref{stabcondmilstein} become
\begin{eqnarray}
x + \frac12 y & < & 0\,, \label{contstabcondxy_m}\\
x + \frac12 y  + \frac12 (1-2\theta) x^2  &<& 0\,,\label{Marstabcondxy_m}\\
x + \frac12 y  + \frac12 \left( (1-2\theta) x^2 + \frac12 \frac1m y^2  + \frac14 y^2 (m-1)^2  \right) &<& 0\,.\label{Milstabcondxy_m}
\end{eqnarray}
Note that for $m=1$ and $m=2$ condition \eqref{Milstabcondxy_m} yields the same stability region.

To interpret the figures below, observe that given a test equation with parameter values $\lambda$ and $\mu_1$ and $m=1$, the point $(x,y) = (\lambda,\mu_1^2)$ corresponds to the choice of step-size $h=1$. Then varying the step-size $h$ corresponds to moving along the ray that connects $(\lambda,\mu_1^2)$ with the origin, where going on this ray in the direction of the origin corresponds to decreasing the step-size $h$. For $m>1$ the scaling is appropriately adapted.

\begin{figure}[htbp]\centering
\includegraphics[clip,bb=-110 110 720 750,width=0.80\textwidth]{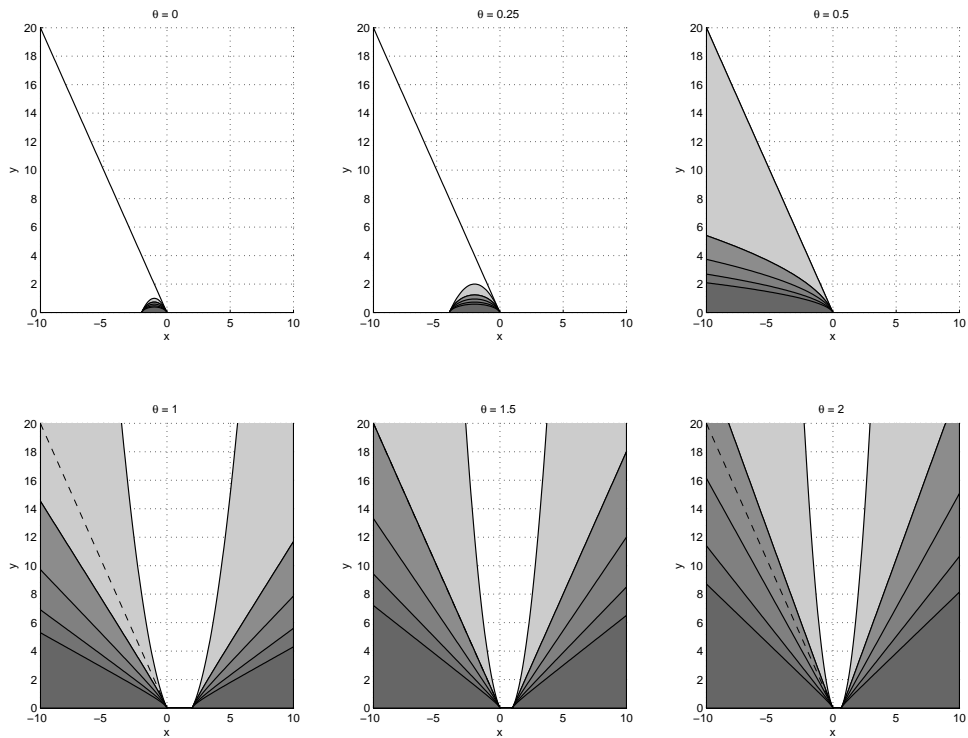}\vspace{-0.4cm}
\caption{\small SDE \eqref{testsde} with $\mu_1 = \mu_2 = \dots = \mu_m$ and $x=h \lambda$ and $y = h m \mu_1^2$: Mean-square stability regions for the zero solutions of the SDE (white area with dashed border), the $\theta$-Maruyama method (light-grey area) and the $\theta$-Milstein method for $m=2,3,4,5$ (light dark-grey area to dark dark-grey areas). }\label{fig.stabregM}
\end{figure}

\medskip%
Figure~\ref{fig.stabregM} shows the mean-square stability regions of the zero solutions of the SDE (white area with a dashed border), the $\theta$-Maruyama approximation (light-grey area) and the $\theta$-Milstein approximation (light dark area to dark-grey area) for different values of the method parameter $\theta$ and illustrates how the stability region of the $\theta$-Milstein method decreases with $m$ and compares to that of the $\theta$-Maruyama method. In particular, one can observe that the mean-square stability regions for the $\theta$-Milstein method are always smaller than the stability region of the $\theta$-Maruyama method. For the $\theta$-Maruyama method the figure illustrates the property of $A$-stability of the method for $\theta\geq \frac12$ for all $m\geq 1$, where as for the $\theta$-Milstein method we can not conclude $A$-stability for a particular $\theta$ and all $m\geq1$.  

\begin{figure}[htbp]\centering
\includegraphics[clip,bb=-110 110 720 740,width=0.80\textwidth]{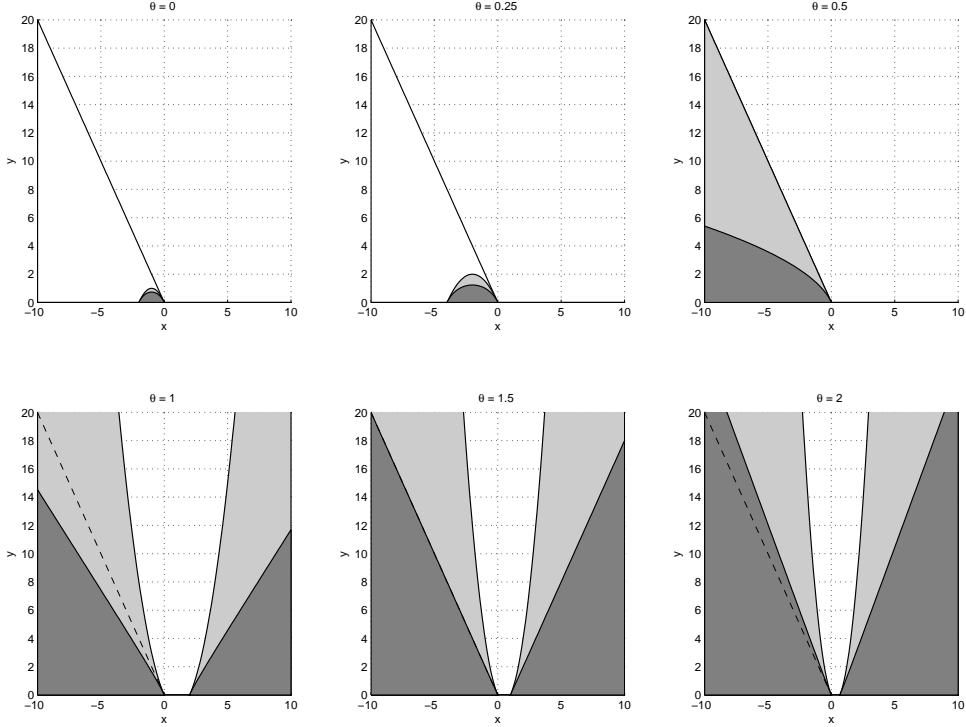}\vspace{-0.4cm}
\caption{\small SDE \eqref{testsde} with $m =1$ and $x=h \lambda$ and $y = h \mu_1^2$: Mean-square stability regions for the zero solutions of the SDE (white area with dashed border), the $\theta$-Maruyama method (light-grey area) and the $\theta$-Milstein method (dark-grey area). }\label{fig.stabregA}
\end{figure}

In \cite[Section 4]{DH00a} and \cite[Section 3]{DH00} the stability regions of the $\theta$-Maruyama method and the $\theta$-Milstein method have been plotted separately for the scalar SDE \eqref{testsde} and real-valued $\lambda, \mu_1$. To emphasise the different stability properties of both methods even in the case of $m=1$, we provide in Figure~\ref{fig.stabregA} a plot of the stability regions of both methods together for the same setting as in \cite{DH00a,DH00} and several values of the method parameter $\theta$. In these cases the stability regions of the $\theta$-Maruyama method and the $\theta$-Milstein method coincide with that of the SDE if $\theta \geq 1/2$ and $\theta \geq 3/2$, respectively. When using the Euler-Maruyama method and the standard Milstein method, that is taking $\theta=0$, the methods that are most often used, then it appears from Figure~\ref{fig.stabregA} that both stability regions are quite small but of similar size. This is already mentioned in \cite[p.114]{ArtAv}. However, when taking mean-square accuracy into account and aiming to reduce numerical costs by using the Milstein method with a larger step-size, one might easily be prevented from doing so by the more restrictive stability condition of the Milstein method.

\section{The effect of introducing partial implicitness in the discretisation of the diffusion term}\label{s.sigma}

A brief look at the deterministic case, that is $\mu_r=0$ ($r=1,\dots,m$) in~\eqref{testsde}, \eqref{thetamethod}, \eqref{stabcondtheta}, etc., reminds us that it is by making the approximation implicit and the introduction of the method parameter $\theta$ that one can control the stability properties of the method by suitably choosing $\theta$. In this sense it would be useful to develop Milstein-type methods that incorporate implicit higher order approximations of the diffusion term to have a method parameter for this purpose.

\medskip%
In general, a straightforward introduction of implicitness into approximations of the diffusion term results in the numerical solution becoming unbounded with positive probability, see \cite[Chap 1.3.4]{MilTret04} and \cite{BuTi01} for methods avoiding this problem. In this article our aim is to highlight the effect that an implicit discretisation of the diffusion term can have on the stability properties of the $\theta$-Milstein method and thus we take advantage of the simple structure of the test equation~\eqref{testsde}. The latter results in the fact that in each term in the second last sum in the $\theta$-Milstein method~\eqref{milsteinmethod} the double Wiener integral $\int_{t_i}^{t_i+h} \int_{t_i}^{s} \rd W_{r}(u)\, \rd W_{r}(s)$ can be replaced by $\frac12 h(\xi_{r,i}^2 - 1)$, for $r=1,\ldots,m$.
Then the second last term in the method~\eqref{milsteinmethod} can be written as $\frac12\, h \mu_r^2 X_i\,\xi_{r,i}^2 - \frac12\, h \mu_r^2 X_i$ for each $r$ and we introduce an implicit approximation
with an additional positive method parameter $\sigma$ only in the latter term, which does not contain the random variable $\xi_{r,i}$. Thus we propose, for each $r$, to use
$$ \frac12\, h  \mu_r^2 \Big( X_i\,
\xi_{r,i}^2 - \big(\sigma X_{i+1} + (1-\sigma) X_i\big)\Big) \quad \text{instead of} \quad \frac12\, h \mu_r^2 \Big( X_i\,\xi_{r,i}^2 -  X_i\Big).$$
We emphasise that this represents a (partial) implicit approximation of the diffusion term in contrast to well-known approaches to implicit approximations of the drift term.

\medskip%
Applied to the scalar linear test-equation \eqref{testsde} the $\theta$-$\sigma$-Milstein method then takes the form
\begin{eqnarray}
\nonumber X_{i+1} &=& X_i + h \big( \theta \lambda X_{i+1}
+ (1-\theta) \lambda X_i \big) + \sqrt{h} \sum_{r+1}^m \mu_r X_i\, \xi_{r,i} \\
 && + \frac12\, h \,\sum_{r=1}^m \mu_r^2 \Big( X_i\,
\xi_{r,i}^2 - \big(\sigma X_{i+1} + (1-\sigma) X_i\big)\Big) + \frac12\, h \!\!\!\sum_{r_1,r_2=1 \atop r_1 \neq r_2 }^m \!\!\mu_{r_1}\mu_{r_2} X_i \xi_{r_1,i} \xi_{r_2,i} \,,\quad \label{sigmamilsteinmethod}
\end{eqnarray}
$i=0,\,1,\ldots\,$, with the random variables $\xi_{r,i}$ defined as in Section~\ref{s.cc}.

\begin{remark}
The choice of $\sigma=\theta$ yields the $\theta$-Milstein approximation of the Stratonovich-SDE $\rd X(t) = (\lambda X(t) - \frac12 \sum_{r=1}^m \mu_r^2 X(t)) \rd t + \sum_{r=1}^m \mu_r X(t) \circ \rd W_r(t)$.
\end{remark}

For the stability analysis we follow the same procedure as in Section~\ref{s.dc} and first rewrite~\eqref{sigmamilsteinmethod} as a one-step recurrence
\begin{equation}\label{sigmarecurrence}
X_{i+1} = \frac1d \bigg( \widehat{a} + \sum_{r=1}^m b_r \,\xi_{r,i} +
\!\!\sum_{r_1,r_2=1}^m \!\! c_{r_1,r_2} \,\xi_{r_1,i} \xi_{r_2,i}\bigg) X_i,
\end{equation}
where the parameters $\widehat{a},b_r,c_{r_1,r_2}$ and $d$ are given by
\begin{eqnarray*}\widehat{a} := 1 + (1-\theta) h \lambda - \frac12 h (1-\sigma) \sum_{r=1}^m \mu_r^2\,,\qquad b_r := \sqrt{h} \mu_r\,, \\
 c_{r_1,r_2} :=\frac12 h \mu_{r_1}\mu_{r_2}\,, \qquad
  d := {1- \theta h \lambda + \frac12 \sigma h \sum_{r=1}^m \mu_r^2}\,.
\end{eqnarray*}
We define $a := \widehat{a} + \sum_{r=1}^m c_{r,r}$ and assume that $(1- \theta h \lambda + \frac12 \sigma h \sum_{r=1}^m \mu_r^2) \neq 0$ to guarantee the existence of a solution to Equation~\eqref{sigmarecurrence}. Now squaring and taking the expectation yields
\begin{eqnarray}
\E |X_{i+1}|^2 &=&
 \frac{1}{|d|^2} \bigg( |a|^2 + \sum_{r=1}^m |b_r|^2 + 2 \sum_{r=1}^m |c_{r,r}|^2 + |c_{sum}|^2\bigg) \, \E |X_i|^2\,,
\label{sigmamilsteinrecurrence}
\end{eqnarray}
where $c_{sum} = \sum_{r_1, r_2=1,\, r_1 \neq r_2 }^m  c_{r_1,r_2}$.

\medskip%
Hence the zero-solution of the stochastic difference equation~\eqref{sigmarecurrence} is asymptotically mean-square stable if and only if the factor  on the right hand side of ~\eqref{sigmamilsteinrecurrence} is less than 1.
Rewriting this condition in terms of $\lambda,\,\mu_r,\,h,\,\theta,\,\sigma$ and rearranging, we obtain the following result:

\begin{lemma}\label{c.stabcondsigma}
The zero solution of the stochastic difference equation given by the $\theta$-$\sigma$-Milstein method~\eqref{sigmamilsteinmethod} applied to the scalar linear test equation~\eqref{testsde} is asymptotically
mean-square stable if and only if
\begin{equation}\label{stabcondsigma}
 \Re(\lambda) + \frac12 \sum_{r=1}^m |\mu_r|^2 + \frac12 h ( 1- 2\theta) |\lambda|^2 + \frac14 h \sum_{r=1}^m|\mu_r|^4 +
 \frac18 h |\mu_{sum}|^2 + \frac12 \sigma h \sum_{r=1}^m \Re(\lambda \mu_r^2)  < 0\,,
\end{equation}
where again $\mu_\text{sum} := \sum_{r_1, r_2=1,\, r_1 \neq r_2 }^m \mu_{r_1} \mu_{r_2}$.
\end{lemma}
We note that the first terms in the left-hand side of~\eqref{stabcondsigma} are equal to
the left-hand side of~\eqref{stabcondmilstein}. The additional term $\frac12 \sigma h \sum_{r=1}^m \Re(\lambda \mu_r^2) $ is negative
for $\Re(\lambda) <0$, i.e., when the stability condition~\eqref{zerostabcond} for the test equation~\eqref{testsde} is satisfied. Thus, the stability condition~\eqref{stabcondsigma} in Lemma~\ref{c.stabcondsigma} is less restrictive than the condition \eqref{stabcondmilstein} for the $\theta$-Milstein method.

\medskip%
Denoting by $S_{\theta\emph{-}\sigma\emph{-Mil}}(\theta,\sigma,h)$ the stability region of the $\theta$-$\sigma$-Milstein method, analogous considerations as for Corollary~\ref{c.Milstein} yield the next result.

\begin{corollary}\label{CorMilstein} Define
\begin{equation*}
\tilde\theta_{\mu_1, \ldots, \mu_m} \,:=\, \frac12 + \sum_{r=1}^{m}
\frac{|\mu_r|^4}{4\,|\lambda|^2} + \frac{|\mu_\text{sum}|^2}{8\, |\lambda|^2} + \frac{\sigma \sum_{r=1}^m \Re(\lambda \mu_r^2)}{2\, |\lambda|^2}\,.
\end{equation*}
For all $h> 0$ it holds that
\[ \begin{array}{rcll}
S_{\theta\emph{-}\sigma\emph{-Mil}}(\theta,\sigma,h)&\subset & S_\text{SDE}& \quad \mbox{for} \quad  0 \leq \theta+\sigma < \tilde\theta_{\mu_1, \ldots, \mu_m}\,,\\
S_{\theta\emph{-}\sigma\emph{-Mil}}(\theta,\sigma,h)&=& S_\text{SDE} & \quad \mbox{for} \quad \theta +\sigma= \tilde\theta_{\mu_1, \ldots, \mu_m}\,,\\
S_{\theta\emph{-}\sigma\emph{-Mil}}(\theta,\sigma,h) &\supset& S_\text{SDE}& \quad
\mbox{for}\quad  \theta+\sigma > \tilde\theta_{\mu_1, \ldots, \mu_m}\,.
\end{array}\]
Then, assuming that $(\lambda,\mu_1, \dots, \mu_m) \in S_\emph{SDE}$, for $0 \leq \theta + \sigma < \tilde\theta_{\mu_1, \ldots, \mu_m}$
the zero solution of the $\theta$-$\sigma$-Milstein method is asymptotically mean-square stable if and only if
\begin{equation}\label{hrestrictionsigma}
 h < \frac{ -2 (\Real(\lambda) + \frac12 \sum_{r=1}^m |\mu_r|^2  + \frac12 \sigma \sum_{r=1}^m \Re(\lambda \mu_r^2) ) }{(1-2\theta)|\lambda|^2 + \frac12 \sum_{r=1}^m |\mu_r|^4 + \frac14 |\mu_\text{sum}|^2} \,.
\end{equation}
In particular, condition \eqref{hrestrictionsigma} then implies that for $\theta> \frac12$ and $ (1-2\theta)|\lambda|^2 + \frac12 \sum_{r=1}^m |\mu_r|^4 + \frac14 |\mu_\text{sum}|^2< 0 $, the stability condition \eqref{stabcondsigma} for the $\theta$-Milstein method is satisfied for any step-size $h > 0$\,. 
\end{corollary}

As in Section~\ref{s.comp} we illustrate the stability regions of the $\theta$-$\sigma$-Milstein method applied to the test SDE \eqref{testsde} with a single noise ($m=1$) in terms of real-valued coefficients $\lambda, \mu_1 \in \R$. In this case we can find a bound on $\tilde\theta_{\mu_1, \ldots, \mu_m}$ as $3/2$. Again, we set $x= h \lambda$ and $y= h \mu_1^2$, and the stability inequality \eqref{stabcondsigma} in terms of $x$ and $y$ reads
\begin{equation*}
   x + \frac12 y  + \frac12 (1-2 \theta) x^2 + \frac14 y^2 +\frac12 \sigma x y <  0\,.
\end{equation*}

\begin{figure}[htbp]\centering\vspace{0.8cm}
\includegraphics[clip,bb=-110 450 720 740,width=0.80\textwidth]{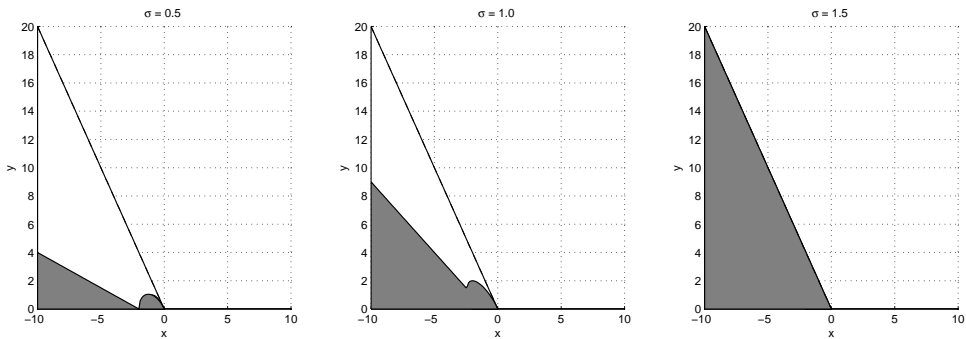}\vspace{-0.4cm}
\caption{\small Mean-square stability regions for the zero solutions of the SDE (white area), the $\theta$-Maruyama method (light-grey area) and the $\theta$-$\sigma$-Milstein method (dark-grey area). Here $\theta = 0$ and $\sigma=0.5, 1, 1.5$.}\label{fig.stabregC}\vspace{0.5cm}
\end{figure}

\begin{figure}[htbp]\centering
\includegraphics[clip,bb=-110 450 720 740,width=0.80\textwidth]{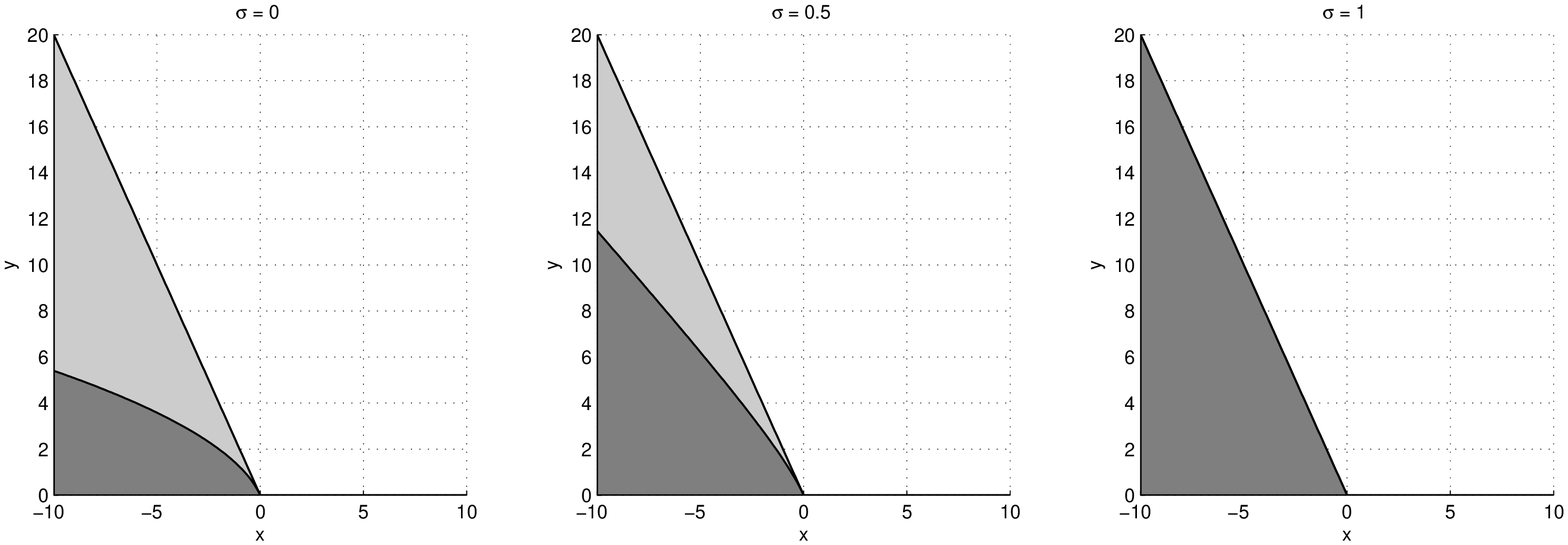}\vspace{-0.4cm}
\caption{\small Mean-square stability regions for the zero solutions of the SDE (white area and dashed line), the $\theta$-Maruyama method (light-grey area) and the $\theta$-$\sigma$-Milstein method (dark-grey area). Here $\theta = 0.5$ and $\sigma=0, 0.5, 1$.}\label{fig.stabregD}
\end{figure}

\begin{figure}[htbp]\centering
\includegraphics[clip,bb=-110 450 720 740,width=0.80\textwidth]{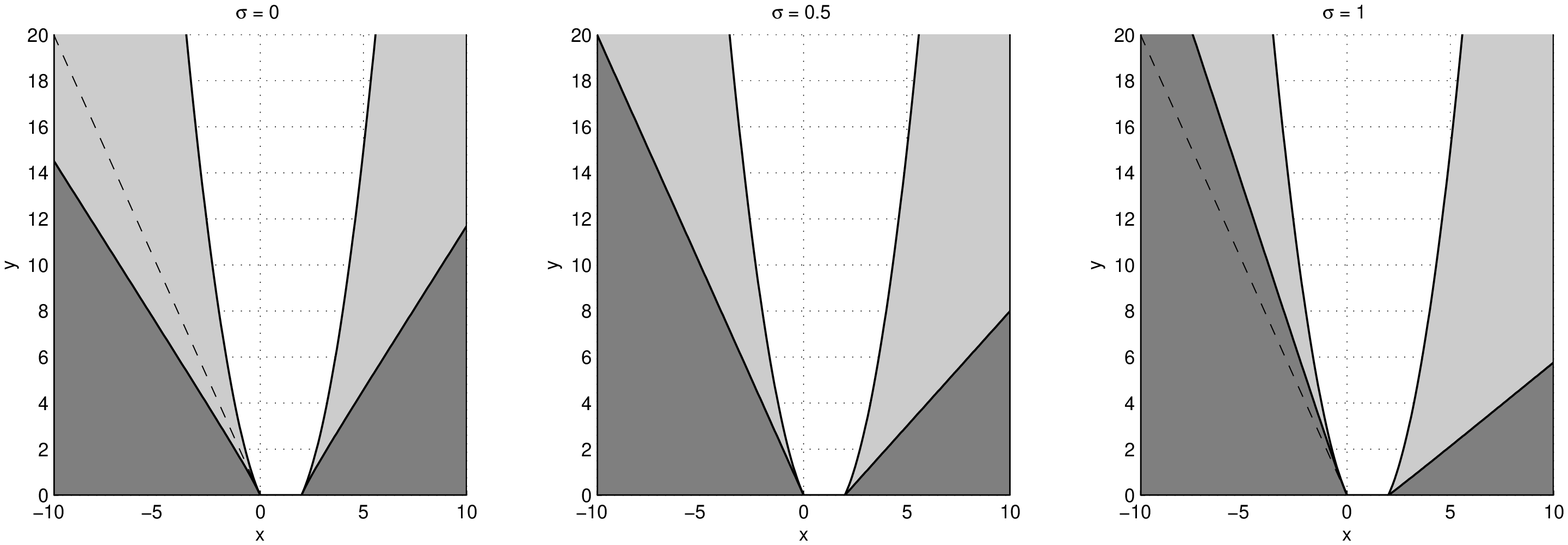}\vspace{-0.4cm}
\caption{\small Mean-square stability regions for the zero solutions of the SDE (dashed line), the $\theta$-Maruyama method (light-grey area) and the $\theta$-$\sigma$-Milstein method (dark-grey area). Here $\theta = 1$ and $\sigma=0, 0.5, 1$.}\label{fig.stabregE}
\end{figure}

\medskip%
Figures \ref{fig.stabregC} - \ref{fig.stabregE} compare the stability regions for the SDE \eqref{testsde} with $m=1$, the $\theta$-Maruyama method and the $\theta$-$\sigma$-Milstein method. In each Figure we have fixed the parameter $\theta$ and only the parameter $\sigma$ varies. For $\sigma = 0$ the plots would correspond to those in Figure~\ref{fig.stabregA}. We can see that the stability region of the $\theta$-$\sigma$-Milstein method (dark-grey area) is becoming larger when increasing the parameter $\sigma$. Moreover, the stability region coincides with the stability region of the test equation  \eqref{testsde} (white area) if $\theta + \sigma \geq 3/2$.

\begin{figure}[htbp]\centering
\includegraphics[clip,bb=-110 110 720 750,width=0.80\textwidth]{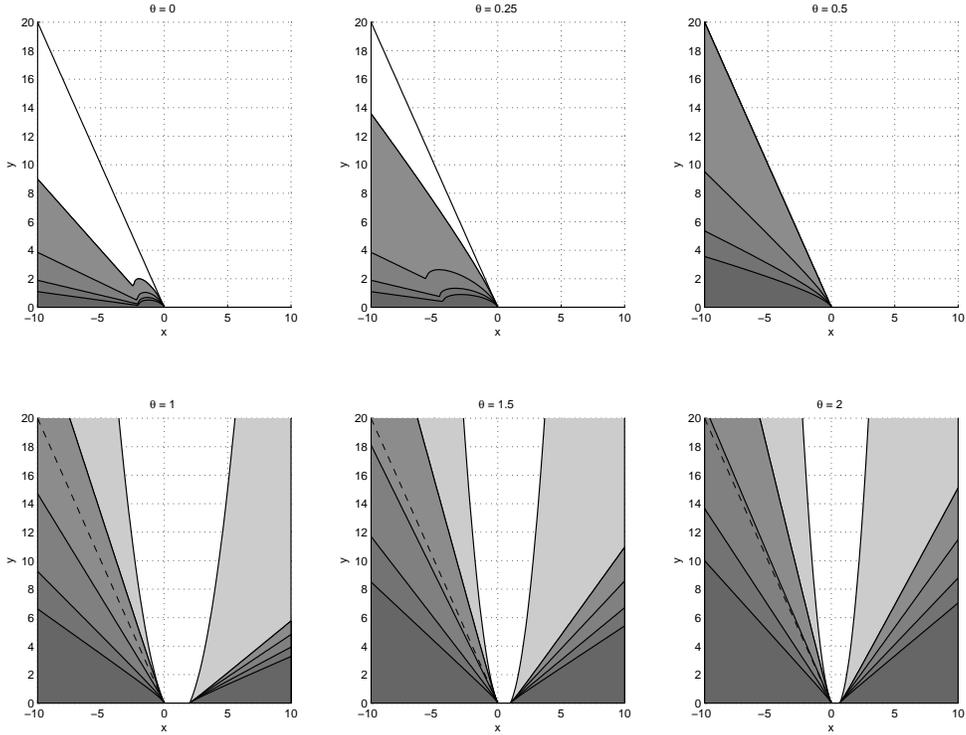}\vspace{-0.4cm}
\caption{\small SDE \eqref{testsde} with $\mu_1 = \mu_2 = \dots = \mu_m$ and $x=h \lambda$ and $y = h m \mu_1^2$: Mean-square stability regions for the zero solutions of the SDE (white area with dashed border), the $\theta$-Maruyama method (light-grey area) and the $\theta-\sigma$-Milstein method for $\sigma = 1$ and $m=2,3,4,5$ (light dark-grey area to dark dark-grey areas). }\label{fig.stabregF}\vspace{-1cm}
\end{figure}

Further, we consider again the SDE \eqref{testsde} with identical noise intensities $\mu_1=\mu_2=\ldots=\mu_m$ and plot in Figure \ref{fig.stabregF} the stability regions for $m=2,\ldots,5$ (again the condition below is the same for $m=1$ and $m=2$) and $\sigma=1$ corresponding to the scaling $x=h\lambda$ and $y=h\,m\,\mu_1^2$, where condition \eqref{stabcondsigma} now reads
\begin{equation*}
   x + \frac12 y  + \frac12 (1-2 \theta) x^2 + \frac14 \frac{1}{m} y^2 + \frac{1}{8} y^2 \,(m-1)^2  +\frac12 \sigma x y <  0\,.
\end{equation*}

Again the stability region of the $\theta$-$\sigma$-Milstein method decreases with growing $m$, as in the case of the $\theta$-Milstein method, which prevents to conclude $A$-stability for a particular value of $\theta + \sigma$ and all $m\geq 1$. However, the introduction of the (partial) implicitness in the diffusion term clearly provides an improvement of the stability behaviour for the $\theta$-$\sigma$-Milstein method.

\section{Numerical experiments}\label{s.num}

In this section we aim to illustrate the effect that the choice of method and the choice of the method parameters $\theta$ and $\sigma$ have on practical simulation runs. We consider the test equation \eqref{testsde} with the following parameters: $m=3$, $\lambda= -2$, $\mu_1=1, \mu_2=-1, \mu_3 = 1$, $X_0 = 0.1$, and integrate for $t\in [0,10]$. The equation has the explicit solution $X(t)= 0.1
\cdot \exp((-2 -\frac{3}{2})t + W_1(t) - W_2(t) + W_3(t))$.

\begin{figure}[htbp]\centering\vspace{0.5cm}
\includegraphics[clip,bb=-110 450 720 740,width=0.90\textwidth]{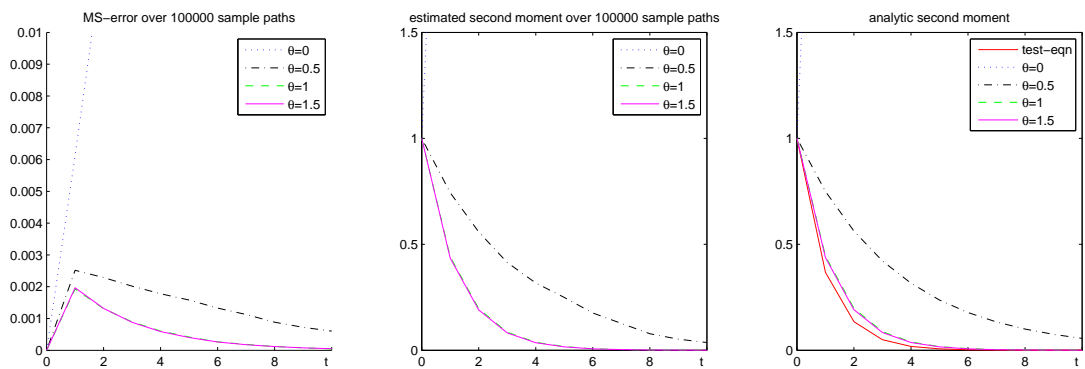}\vspace{-0.4cm}
\caption{\small MS-error, estimated and analytic second moments using the $\theta$-Maruyama method with step-size $h=1$.}\label{fig.stabresA}\vspace{0.5cm}
\end{figure}

\begin{figure}[htbp]\centering
\includegraphics[clip,bb=-110 450 720 740,width=0.90\textwidth]{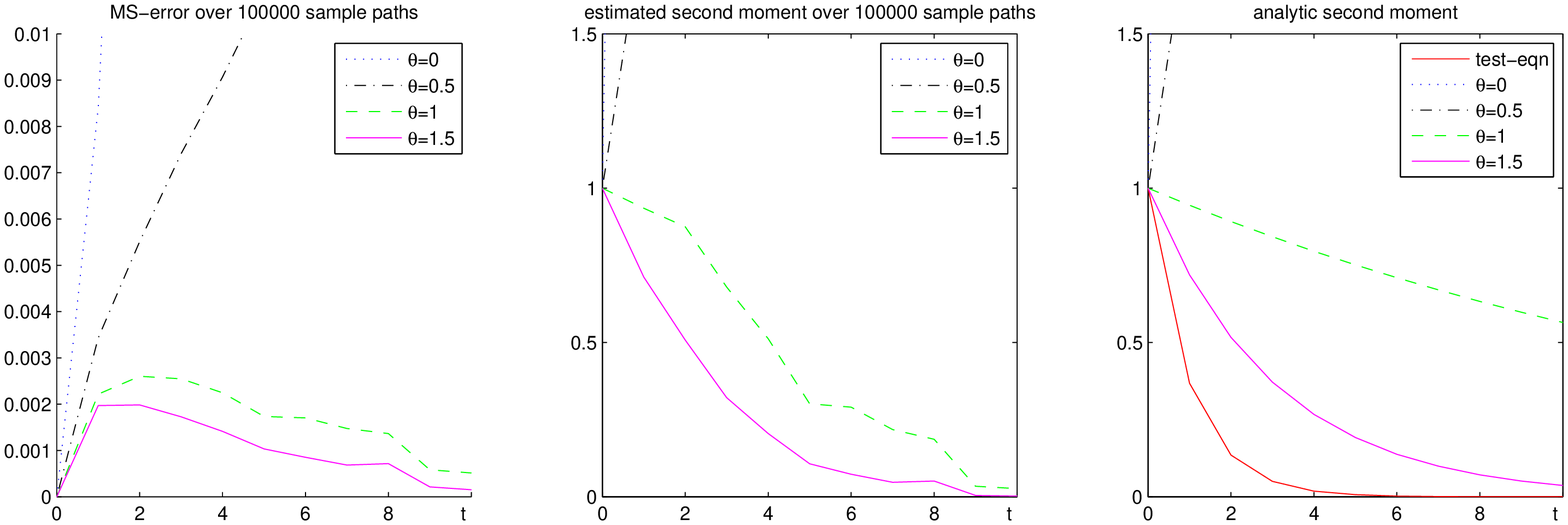}\vspace{-0.4cm}
\caption{\small MS-error, estimated and analytic second moments using the $\theta$-Milstein method with step-size $h=1$.}\label{fig.stabresC}
\end{figure}

\begin{figure}[htbp]\centering
\includegraphics[clip,bb=-110 450 720 740,width=0.90\textwidth]{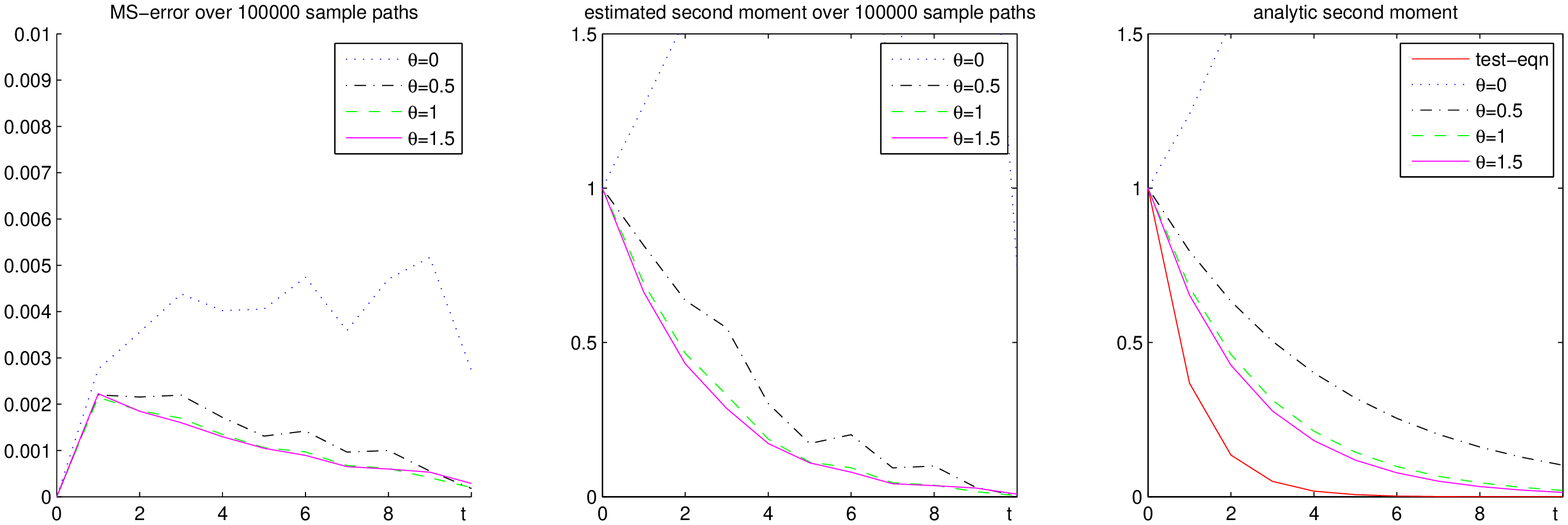}\vspace{-0.4cm}
\caption{\small MS-error, estimated and analytic second moments using the $\theta$-$\sigma$-Milstein method with step-size $h=1$ and $\sigma =1$.}\label{fig.stabresD}
\end{figure}

\begin{figure}[htbp]\centering
\includegraphics[clip,bb=-110 450 720 740,width=0.90\textwidth]{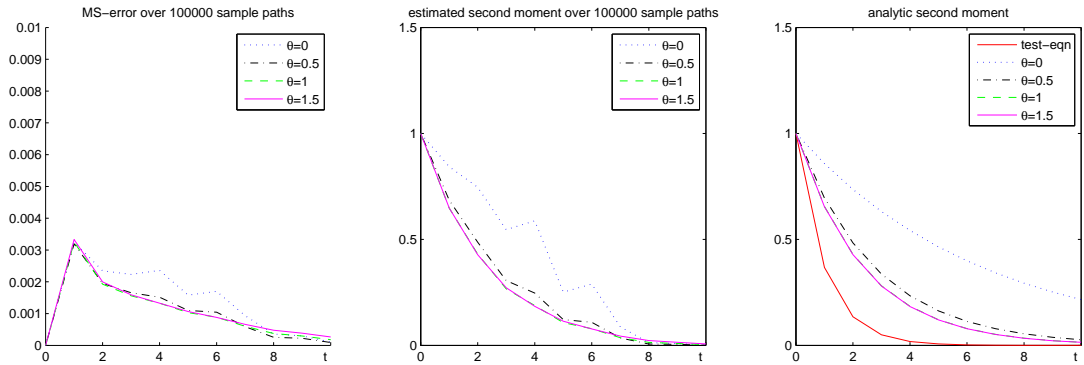}\vspace{-0.4cm}
\caption{\small MS-error, estimated and analytic second moments using the $\theta$-$\sigma$-Milstein method with step-size $h=1$ and $\sigma =1.5$.}\label{fig.stabresE}
\end{figure}

\medskip%
The following four figures show numerical simulation studies performed with the $\theta$-Maruyama method (Figure~\ref{fig.stabresA}), the $\theta$-Milstein method (Figure~\ref{fig.stabresC}), the $\theta$-$\sigma$-Milstein method with $\sigma =1$ (Figure~\ref{fig.stabresD}) and the $\theta$-$\sigma$-Milstein method with $\sigma =1.5$ (Figure~\ref{fig.stabresE}), all simulations done with a fixed step-size $h=1$ over the interval $[0,10]$ with gridpoints $t_i=i\,h,\,i=0,\ldots,10$. The parameter $\theta$ varies as $\theta = 0, \, 0.5,\, 1,\, 1.5$.

\begin{remark}
The analytical results in the previous sections also suggest choices of step-sizes, such that the zero solution of any one of the numerical methods for a fixed set of parameters is asymptotically mean-square stable and obviously it is possible to illustrate this by performing numerical experiments with fixed parameter sets and varying the step-size $h$. However, the focus of this article is rather on comparing qualitative properties of \emph{methods} than on how individual methods behave for various step-sizes.
\end{remark}

Each of the four figures below consists of three pictures.
The underlying idea of a stability analysis of numerical methods is to provide guidance for choosing method parameters and a step-size such that the numerical solution represents a good approximation of the true solution. (Note that convergence of a method only guarantees this in the limit for $h\rightarrow 0$.) Thus, we plot in the left picture of each figure the mean-square error $e$ between the exact solution $X(t_i)$ and numerical solution $X_i$ for the corresponding method and choice of method parameters.
This allows to compare the impact of the choice of method on actually computing solutions of SDEs. However, the quantity of interest in the stability definitions \ref{d.stabil} and \ref{d.stabildiff} are the second moments of the analytical and numerical solutions. Therefore the other two pictures present plots of the second moments of the solutions of \eqref{testsde}, \eqref{thetamethod}, \eqref{milsteinmethod} and \eqref{sigmamilsteinmethod}, in the middle picture estimated from the numerical simulations of  \eqref{thetamethod}, \eqref{milsteinmethod} and \eqref{sigmamilsteinmethod}, in the right picture computed form the analytical solutions. This allows to compare the effect of the choice of method and parameters on the behaviour of the second moments analytically and also, how well this is approximated by the numerical methods.

\medskip%
In detail, the figures present the following quantities, where $(X(t_i,\omega_j))_{i=0,\dots,N}$ and $(X_i(\omega_j))_{i=1,\dots,N}$ denote the values on grid points of a trajectory of the explicit solution of \eqref{testsde} with the above parameters and of the numerical trajectory produced by one of the methods, respectively:
\begin{eqnarray}
\text{MS-error}(t_i) & = & \Bigg(\frac1M \sum_{j=1}^M (X(t_i,\omega_j) - X_i(\omega_j))^2 \Bigg)^{1/2} \approx (\E |X(t_i) - X_i|^2 )^{1/2}\,, \nonumber \\
\widehat{\E}(X_i^2) & = & \frac1M \sum_{j=1}^M X_i^2(\omega_j) \quad \text{ for all } i \,, \label{estsecmom}\\
\E(X^2(t_i)) & = & X_0^2 \cdot \exp\left\{ (2 \lambda + \sum_{r=1}^m \mu_r^2) \cdot t_i \right\} \quad \text{ for all } i \,, \label{expsecmom}\\[1em]
\E(X_i^2) & = & s^i \cdot X_0^2 \quad \text{ for all } i \,. \label{numsecmom}
\end{eqnarray}
The expression \eqref{estsecmom} for $\widehat{\E}(X_i^2)$ represents an estimator for the second moment of the numerical approximation, the expression \eqref{expsecmom} for $\E(X^2(t_i))$ is the solution of the deterministic ODE $\E(X^2(t))' = (2 \lambda + \sum_{r=1}^m \mu_r^2)\, \E(X^2(t))$ (see proof of Lemma~\ref{l.contprob}) at discrete time-points and the last expression for $\E(X_i^2)$ is the result of applying the recurrences \eqref{sqthetarecurrence}, \eqref{milsteinrecur} and \eqref{sigmamilsteinrecurrence} where $s$ denotes the factor in front of $\E(X_i^2)$ in the right-hand side of the corresponding recurrence equation. The number of trajectories computed for the above quantities is $M=100000$ for each simulation.

\medskip%
Figures~\ref{fig.stabresA} to~\ref{fig.stabresE} illustrate the behaviour of the methods for the above set of parameters, in particular the same fixed step-size $h=1$ and different choices of $\theta$ and $\sigma$. The $\theta$-Maruyama method provides reliable approximations for $\theta \geq 0.5$, but not for $\theta=0$, the $\theta$-Milstein method in same setting \emph{does not} provide reliable approximations for  $\theta=0$ and $\theta=0.5$. In fact, the numerical solutions for $\theta=0.5$ in this setting diverge. We can observe the improvement in the stability behaviour when using the $\theta$-$\sigma$-Milstein method introduced in the previous section. In contrast to the $\theta$-Milstein method this method produces reliable approximations for the choice of $\theta= 0.5$ when also setting $\sigma=1$. Further, for the choice of $\sigma=1.5$ the $\theta$-$\sigma$-Milstein approximation behaves satisfactorily for all values of $\theta$, thus the implicit term in the diffusion approximation is effectively stabilising the method.


\section{Conclusions}\label{s.conclusions}

A linear stability analysis has been performed for the $\theta$-Maruyama method and the $\theta$-Milstein method, using a linear test equation with several multiplicative noise terms. We have obtained stability conditions guaranteeing asymptotic mean-square
stability of the zero solution of the stochastic difference equations resulting from both types of methods applied to the test SDE. Comparing the stability conditions \eqref{stabcondtheta} and \eqref{stabcondmilstein} for the $\theta$-Maruyama and the $\theta$-Milstein method, respectively, it is quite obvious that the latter is more restrictive than the former, due to the fourth and fifth term in \eqref{stabcondmilstein}, i.e.~$\frac14 h \sum_{r=1}^m |\mu_r|^4 + \frac18 h |\mu_{sum}|^2$, which is always non-negative. In particular, this term is the more restrictive, the larger the parameters in the diffusion term in~\eqref{testsde} are. Now, in the case that these parameters are small, that is the so-called small noise case, it is well known from the corresponding mean-square convergence analysis that Maruyama-type methods provide sufficiently accurate methods for practicable choices of step-sizes, see~\cite{BW06,MilTret,Si08}. Milstein-type methods include higher order approximations of the diffusion term with the aim that they are accurate methods for more practicable choices of step-sizes just in the case that the diffusion term is large. In other words, when dealing with SDEs with larger diffusion terms, numerical efficiency considerations would suggest using a Milstein-type method with a larger step-size rather than a Maruyama-type method with a small step-size to obtain numerical approximations of the solution of an SDE with a similar accuracy. However, Corollary~\ref{CorMilStab} indicates that in this case there may be a trade-off and one is faced with restrictions on the step-size for the Milstein-type method due to stability reasons. Further, we have shown that the precise stability region of the $\theta$-Milstein method depends on the number and magnitude of the noise terms, whereas the stability region of the $\theta$-Maruyama method is independent of them. In particular, it is not possible to define a value $\theta_{\text{bound}}$ such that the $\theta$-Milstein method can be called $A$-stable for all $\theta>\theta_{\text{bound}}$ for the class of SDEs \eqref{testsde} for an arbitrary number of driving Wiener processes. For the $\theta$-Maruyama method $\theta_{\text{bound}}=\frac12$, as in the deterministic case and the case of the test equation \eqref{testsde} with $m=1$. We provide a modified Milstein-type method with a partially implicit diffusion approximation and demonstrate that the resulting stability behaviour can be controlled more favourably. The results highlight that it is necessary to include multi-dimensional noise into test equations and to study their effects on the practical behaviour of the methods.

\section*{Acknowledgement}\label{s.thanks}

We thank Martin Riedler (Heriot-Watt University, Edinburgh) and Lukasz Szpruch (University of Strathclyde, Glasgow) for fruitful discussions.

\end{document}